\documentclass[12pt]{amsart}
\usepackage{amssymb}
\usepackage{amsfonts}
\usepackage{latexsym}
\usepackage{amscd}

\newfont{\gd}{eufm10 scaled \magstep1}
\newfont{\gs}{eufm7 scaled \magstep1}
\newfont{\gss}{eufm5 scaled \magstep1}

\newcommand{\length}{\textrm{len}}

\newcommand{\A}{{\mathbb{A}}}
\newcommand{\N}{{\mathbb{N}}}
\newcommand{\Z}{{\mathbb{Z}}}

\newcommand{\ndr}{\mbox{\rm Endr}}
\newcommand{\id}{\mbox{\rm id}}
\newcommand{\D}{\mathcal{D}}
\newcommand{\Dtil}{\widetilde{D}}

\newcommand{\ttil}{\tilde{t}}
\newcommand{\sop}{S^{\text{\rm op}}}
\newcommand{\sinvt}{S^{-1}T}
\newcommand{\sinvs}{S^{-1}S}
\newcommand{\amod}{A\mbox{\rm-Mod}}
\newcommand{\ndz}{\mbox{\rm End}_{\mathbb{Z}}}
\newcommand{\sopat}{\sop*_\alpha A*_\alpha T}
\newcommand{\sopas}{\sop*_\alpha A*_\alpha S}
\newcommand{\Zplus}{\Z^+}
\newcommand{\zpazp}{\Zplus*_\alpha A*_\alpha \Zplus}
\newcommand{\zpbzp}{\Zplus*_\alpha B*_\alpha \Zplus}
\newcommand{\atpmalpha}{A[t_+,t_-;\alpha]}
\newcommand{\btpmalpha}{B[t_+,t_-;\alpha]}
\newcommand{\dztpmalpha}{D_0[t_+,t_-;\alpha]}
\newcommand{\aut}{\mbox{\rm Aut}}
\newcommand{\alsom}{{\widehat {\alpha }}}
\newcommand{\sopasprime}{\sop*_{\alpha'}(eAe)*_{\alpha'}S}
\newtheorem{lemma}{Lemma}[section]
\newtheorem{corollary}[lemma]{Corollary}
\newtheorem{theorem}[lemma]{Theorem}
\newtheorem{proposition}[lemma]{Proposition}

\newtheorem{definition}[lemma]{Definition}
\newtheorem{example}[lemma]{Example}

\newtheorem{noname}[lemma]{}

\begin{document}

\title[Fractional skew monoid rings]{Fractional skew
monoid rings}%
\author{P. Ara}
\address{Departament de Matem\`atiques, Universitat Aut\`onoma de Barcelona, 08193 Bellaterra (Barcelona), Spain.}
\email{para@mat.uab.es}
\author{M.A. Gonz\'alez-Barroso}
\address{Departamento de Matem\'aticas, Universidad de C\'adiz,
Apartado 40, 11510 Puerto Real (C\'adiz),
Spain.}\email{mariangeles.gonzalezbarroso@alum.uca.es}
\author{K.R. Goodearl}
\address{Department of Mathematics, University of California Santa
Barbara CA 93106, U.S.A.}
\email{goodearl@math.ucsb.edu}
\author{E. Pardo}
\address{Departamento de Matem\'aticas, Universidad de C\'adiz,
Apartado 40, 11510 Puerto Real (C\'adiz),
Spain.}\email{enrique.pardo@uca.es}

\thanks{The first author was partially supported by the DGI
and European Regional Development Fund, jointly, through Project
BFM2002-01390, the second by an FPU fellowship of the Junta de
Andaluc\'{\i}a, the third by an NSF grant, and the fourth by DGESIC
grant BFM2001--3141. Also, the first and fourth authors are
partially supported by the Comissionat per Universitats i Recerca
de la Generalitat de Catalunya, and the second and fourth by PAI
III grant FQM-298 of the Junta de Andaluc\'{\i}a.}



\begin{abstract}
Given an action $\alpha$ of a monoid $T$ on a ring $A$ by ring
endomorphisms, and an Ore subset $S$ of $T$, a general
construction of a fractional skew monoid ring $\sop*_\alpha
A*_\alpha T$ is given, extending the usual constructions of skew
group rings and of skew semigroup rings. In case $S$ is a
subsemigroup of a group $G$ such that $G=S^{-1}S$, we obtain a
$G$-graded ring $\sopas$ with the property that,
for each $s\in S$, the $s$-component contains a left invertible
element and the $s^{-1}$-component contains a right invertible
element. In the most basic case, where $G=\Z$ and $S=T=\Zplus$,
the construction is fully determined by a single ring endomorphism
$\alpha$ of $A$. If $\alpha $ is an isomorphism onto a proper
corner $pAp$, we obtain an analogue of the usual skew Laurent
polynomial ring, denoted by $\atpmalpha$. Examples of this
construction are given, and it is proven that several classes of
known algebras, including the Leavitt algebras of type $(1,n)$,
can be presented in the form $\atpmalpha$. Finally, mild and
reasonably natural conditions are obtained under which $\sopas$ is
a purely infinite simple ring.
\end{abstract}
\maketitle

\section*{Introduction}

Let $\alpha\colon G\rightarrow \aut (A)$ be an action of a group
$G$ on a unital ring $A$. A useful construction attached to these
data is the {\it skew group ring} $A*_{\alpha}G$, see \cite{Mon}
and \cite{Passman}. This is the ring of formal expressions $\sum
_{g\in G} a_gg$, where $a_g\in A$ and almost all the coefficients
$a_g$ are $0$. Addition is defined componentwise and
multiplication is defined according to the rule $(ag)(bh)=(a\alpha
_g(b))(gh)$. The skew group ring $A*_{\alpha}G$ can also be
defined as the unital ring $R$ such that there are a unital ring
homomorphism $\phi \colon A\rightarrow R$ and a unital monoid
homomorphism $i\colon G\rightarrow R$ from $G$ to the
multiplicative structure of $R$, universal with respect to the
property that $i(g)\phi(a)=\phi(\alpha _g(a))i(g)$ for all $a\in
A$ and all $g\in G$. In his pioneering paper \cite{Paschke},
Paschke gave a construction of a $C^*$-algebraic crossed product
$A\rtimes _{\alpha}\N$ associated to a not necessarily unital
$C^*$-algebra endomorphism $\alpha$ on a $C^*$-algebra $A$.
Paschke's $C^*$-algebraic construction has been generalized to
other semigroups, see \cite{LacaRae}, \cite{Larsenirish},
\cite{Larsencanad} and \cite{LarsenRae}. Moreover, R\o rdam
\cite{rordamclassif} used Paschke's construction together with the
Pimsner-Voiculescu exact sequence associated to an automorphism
\cite[Theorem 10.2.1]{Black} to realize any pair of countable
abelian groups $(G_0,G_1)$ as $(K_0(B),K_1(B))$ for a certain
purely infinite, simple, nuclear separable $C^*$-algebra $B$.

In this paper, we develop a systematic purely algebraic theory of
fractional skew monoid rings with respect to monoid actions on
rings by not necessarily unital ring endomorphisms, in which an
Ore submonoid is inverted. (Recall that a monoid is a semigroup
with a neutral element.) More precisely, we assume the following
data are given (see (\ref{data1}) for the detailed definitions of
the properties):
\begin{enumerate}
\item A monoid $T$ acting on a unital ring $A$ by endomorphisms;
\item A submonoid $S$ of $T$ satisfying the left denominator
conditions, and such that $S$ is left saturated in $T$.
\end{enumerate}
Then a fractional skew monoid ring $\sopat$ is constructed, with
suitable maps from $A$, $\sop$ and $T$ to $\sopat$, which satisfy
a universal property analogous to the one for the skew group ring
described above, see Definition \ref{newring}. It is not difficult
to show that such a ring exists by using a construction with
generators and relations, but it is rather non-obvious to
determine the algebraic structure of $\sopat$. The ring $\sopat$
is best understood by means of its $S^{-1}T$-graded structure,
obtained in Proposition \ref{gradedring}. The structure is
completely pinned down in (\ref{alphasinject}) in the case where
$T$ acts by injective endomorphisms.

The general construction of $\sopat$ is given in Section 1. In the
other sections, we specialize the construction to the case of a
submonoid $S$ of a group $G$ such that $G=S^{-1}S$ (taking
$T=S$), and to an action $\alpha$ of $S$ on $A$ by
corner isomorphisms, meaning that $\alpha _s$ is an isomorphism
from $A$ onto the corner ring $\alpha _s(1)A\alpha _s(1)$ for all
$s\in S$. Several examples of interest are considered in Section 2
in the case where $S=T=\Zplus$. In particular, the Leavitt algebras
$V_{1,n}(k)$ and $U_{1,n}(k)$, already considered by Leavitt,
Skornyakov, Cohn, Bergman and others, are seen here to be
particular cases of our construction.

For $S=T=\Zplus$, the construction is determined by a single
corner isomorphism $\alpha$, and the elements of the fractional
skew monoid ring $R= \zpazp$ can all be written as `polynomials' of
the form
$$r=a_{n}t_+^n+\ldots
+a_{1}t_++a_0+t_{-}a_{-1}+\ldots t_{-}^m a_{-m},$$
 with coefficients $a_i\in A$. Because of this similarity of $R$ with
a skew-Laurent polynomial ring, we shall use the notation $R=
\atpmalpha$.  Using this construction and the
Bass-Heller-Swan-Farrell-Hsiang-Siebenmann Theorem,  the $K_1$
group of these algebras is computed in \cite{AB}.

A general source of interesting examples is provided in Section 3.
Namely, assume that $G$ is a group acting on a ring $A$ by
automorphisms, and that there are a submonoid $S$ of $G$ such that
$G=S^{-1}S$ and a non-trivial idempotent $e$ in $A$ such that
$\alpha _s(e)\le e$ for all $s\in S$. Then the corner ring
$e(A\ast _{\alpha} G)e$ of the skew group ring $A\ast _{\alpha} G$
is isomorphic as a $G$-graded ring to a fractional skew monoid
ring $\sopasprime$ (Proposition
\ref{isomap}). In case $G$ is abelian, all rings $\sopas $ appear
in this way (Proposition \ref{genunit}).

Sections 4 and 5 deal with actions on simple rings. Using a
suitable definition of outer action of a monoid $S$ on a ring $A$,
we prove in Theorem \ref{simplicity} that $\sopas $ is a simple
ring for any outer action $\alpha$ of $S$ on a simple ring $A$.
This is a generalization of a well-known sufficient condition for
simplicity of skew group rings, see \cite[Theorem 2.3]{Mon}.
Section 5 shows that, under mild conditions on $A$ and on the
outer action $\alpha$ of $S$ on $A$, the fractional skew monoid
ring $\sopas $ is a purely infinite simple ring (Theorem
\ref{spi}). In particular, this holds whenever $A$ is either a
simple ultramatricial algebra over some field or a purely infinite
simple ring. The class of purely infinite simple rings has been
recently studied by the first, third and fourth authors in
\cite{agp}, and constitute an important and large class of
relatively well-behaved simple rings. They can be thought of as
the nice rings in the wild universe of the directly infinite
simple rings; see specially \cite[Corollary 2.2 and Theorem
2.3]{agp} for the good behaviour of $K$-theory of purely infinite
simple rings. A further nice property of them has been recently
established by the first author in \cite{exchange}: Every purely
infinite simple ring satisfies the exchange property.

All rings and modules in this paper will be assumed to be unital
unless explicitly noted. (The main exception is the ring $\A$
constructed in Section \ref{skewmonfromskewgrp}.) However, many of
the subrings we deal with will have units different from the unit
of the larger ring; specifically, we will deal with many
\emph{corners} $pAp$ in a ring $A$, where $p$ is an idempotent.
Note that any ring endomorphism $\varepsilon$ of
$A$, even if not unital when considered as a map $A\rightarrow A$,
\underbar{is} unital when viewed as a ring homomorphism
$A\rightarrow \varepsilon(1)A\varepsilon(1)$.

\section{The general construction}

We present the construction of a fractional skew monoid ring in
full generality in this section, and establish the precise graded
structure of this ring. The basic data consist of a ring $A$, a
monoid $T$ acting on $A$ by ring endomorphisms, and a left
denominator set $S\subseteq T$; the fractional skew monoid ring we
construct is graded by $\sinvt$, and its identity component is the
quotient of $A$ modulo the union of the kernels of the
endomorphisms by which $S$ acts.

\begin{noname}\label{data1} {\rm
 We begin by fixing the basic data
needed for our construction; these data and conventions will
remain in force throughout the paper.
 Let $A$ be a (unital) ring, and $\ndr(A)$
the monoid of non-unital (i.e., not necessarily unital) ring
endomorphisms of $A$.

Let $T$ be a monoid and $\alpha: T\rightarrow \ndr(A)$ a monoid
homomorphism, written $t\mapsto \alpha_t$. In general, we will
write $T$ multiplicatively, with its identity element denoted $1$,
but in some applications it will be convenient to switch to
additive notation for $T$. For $t\in T$, set $p_t= \alpha_t(1)$,
an idempotent in $A$. Then $\alpha_t$ can be viewed as a unital
ring homomorphism from $A$ to the corner $p_tAp_t$. For $s,t\in
T$, we have $p_{st}= \alpha_{st}(1)= \alpha_s \alpha_t(1)=
\alpha_s(p_t)$.

Let $S\subseteq T$ be a submonoid satisfying the left denominator
conditions, i.e., the left Ore condition and the monoid version of
left reversibility: whenever $t,u\in T$ with $ts=us$ for some
$s\in S$, there exists $s'\in S$ such that $s't=s'u$. Then there
exists a monoid of fractions, $S^{-1}T$, with the usual properties
(e.g., see \cite[\S1.10]{ClPr} or \cite[\S0.8]{Cohn}).

We shall also assume that $S$ is \emph{left saturated} in $T$:
whenever $s\in S$ and $t\in T$ such that $ts\in S$, we must have
$t\in S$. This assumption means that equality in $\sinvt$ can be
described as follows: if $s_1^{-1}t_1= s_2^{-1}t_2$ for some
$s_i\in S$ and $t_i\in T$, there exist $u_1,u_2\in S$ such that
$u_1s_1= u_2s_2$ and $u_1t_1= u_2t_2$. (The usual denominator
conditions only yield the latter equations for, say, some $u_1\in
S$ and $u_2\in T$. But then $u_2s_2= u_1s_1\in S$, and left
saturation implies $u_2\in S$.) }
\end{noname}

\begin{definition}\label{newring} {\rm The label $\sopat$
stands for a (unital) ring $R$ equipped with a (unital) ring
homomorphism $\phi: A\rightarrow R$ and monoid homomorphisms
$s\mapsto s_-$ from $\sop\rightarrow R$ and $t\mapsto t_+$ from
$T\rightarrow R$, universal with respect to the following
relations:

\begin{enumerate}

\item $t_+\phi(a)= \phi\alpha_t(a)t_+$ for all $a\in A$ and $t\in
T$;

\item $\phi(a)s_-= s_-\phi\alpha_s(a)$ for all $a\in A$ and $s\in
S$;

\item $s_-s_+= 1$ for all $s\in S$;

\item $s_+s_-= \phi(p_s)$ for all $s\in S$.

\end{enumerate}

Note that condition (2) follows from the
others. Given $a\in A$ and $s\in S$, we have $s_+\phi(a)=
\phi\alpha_s(a)s_+$ by (1), and on multiplying each term of this
equation on the left and on the right by $s_-$, we obtain
$\phi(a)s_-= s_-\phi\alpha_s(a)\phi(p_s)= s_-\phi(\alpha_s(a)p_s)$
from (3) and (4), whence (2) follows because $\alpha_s(a)p_s=
\alpha_s(a)$.  }
\end{definition}

\begin{noname}
\label{discussion}{\rm  At this point, we sketch the
\emph{existence} of the ring $\sopat$. The existence of a
ring satisfying the universal property of Definition
\ref{newring} follows from a construction with generators and
relations, which does not use at all any property of $S$; in fact,
$S$ can be an arbitrary subset of $T$. Take $B=A*\Z \langle
t_+,s_-\mid t\in T,\, s\in S\rangle$ to be the free product of $A$
and the free ring on the disjoint union $T\sqcup S$, and let
$i_1\colon A\rightarrow B$ and
$i_2\colon \Z \langle t_+,s_-\mid t\in T,\, s\in
S\rangle\rightarrow B$ be the canonical maps. Let $I$ be the
two-sided ideal of $B$ generated by

\quad(a) $i_2(t_+)i_1(a)- i_1(\alpha_t(a))i_2(t_+)$ for all $a\in
A$ and $t\in T$;

\quad(b) $i_2((tt')_+)-i_2(t_+)i_2(t'_+)$ for all $t,t'\in T$;

\quad(c) $i_2(s_-)i_2(s_+)- i_1(1)$ for all $s\in S$;

\quad(d) $i_2(s_+)i_2(s_-)- i_1(p_s)$ for all $s\in S$.

\noindent Then $\sopat=B/I$ is the ring we are looking for, and
$\phi$ is the composite map $\pi \circ i_1$, where $\pi\colon
B\rightarrow B/I$ is the canonical projection. Note that the
relations
$(ss')_-=s'_-s_-$, for all $s,s'\in S$ such that $ss'\in S$, hold
automatically from (a)--(d) above. Also, we have already observed
that condition (2) in \ref{newring} follows from conditions
(1),(3) and (4), and so it follows from (a)--(d) too.

Rather than introduce a notation for the product in $\sop$, we
view the map $(-)_-$ as a monoid anti-homomorphism $S\rightarrow
R$, so that $(su)_-= u_-s_-$ for $s,u\in S$.

The construction above will also be applied when $A$ is an algebra
over a field $k$ and the ring endomorphisms $\alpha_t$ for $t\in
T$ are $k$-linear. In this case, it is easily checked that $\phi$
maps $k=k\cdot1$ into the center of $B$ (use relations (1),(2)
above and part (c) of the following lemma to see that $\phi(k)$
commutes with each $s_-$ and $t_+$), so that $B$ becomes a
$k$-algebra and $\phi$ becomes a $k$-algebra homomorphism. The
universal property of $B$ then holds also in the category of
$k$-algebras.}
\end{noname}

The following lemma and subsequent results pin down the structure
of $R$. This structure simplifies considerably when the maps
$\alpha_s$ are injective -- see (\ref{alphasinject}).

\begin{lemma}\label{relations1} Let $a,b\in A$, $s,u\in S$, and $t,v\in T$.

{\rm (a)} $s_+\phi(a)s_-= \phi\alpha_s(a)$.

{\rm (b)} $s_-\phi\alpha_s(a)s_+= \phi(a)$.

{\rm (c)} $s_-= s_-\phi(p_s)$ and $t_+= \phi(p_t)t_+$.

{\rm (d)} $s_-\phi(a)t_+= s_-\phi(p_sap_t)t_+$.

{\rm (e)} $s_-\phi(a)t_+= (us)_-\phi\alpha_u(a)(ut)_+$.

{\rm (f)} There exist $x\in S$ and $y\in T$ such that $xt=yu$. For
any such $x,y$,
$$\bigl[ s_-\phi(a)t_+ \bigr] \bigl[ u_-\phi(b)v_+
\bigr]= (xs)_-\phi\bigl(\alpha_x(ap_t)\alpha_y(b)\bigr)(yv)_+.$$
 In particular, $t_+u_-= x_-p_{xt}y_+$.
\end{lemma}

\noindent\emph{Proof}. (a) $s_+\phi(a)s_-= \phi\alpha_s(a)s_+s_-=
\phi\alpha_s(a)\phi(p_s)= \phi(\alpha_s(a)p_s)= \phi\alpha_s(a)$.

(b) This follows from (a) because $s_-s_+= 1$.

(c) $s_-= \phi(1)s_-= s_-\phi\alpha_s(1)= s_-\phi(p_s)$.
Similarly, $t_+= t_+\phi(1)= \phi\alpha_t(1)t_+= \phi(p_t)t_+$.

(d) This is clear from (c).

(e) From (b), we have $\phi(a)= u_-\phi\alpha_u(a)u_+$, and the
desired equation follows because $s_-u_-= (us)_-$.

(f) Note that $(xt)_+(yu)_-= (xt)_+(xt)_-= \phi(p_{xt})=
\phi\alpha_x(p_t)$. Using (e), we get
 \begin{align}
 \nonumber \bigl[ s_-\phi(a)t_+ \bigr] \bigl[ u_-\phi(b)v_+
\bigr]  &= \bigl[ (xs)_-\phi\alpha_x(a)(xt)_+ \bigr] \bigl[
(yu)_-\phi\alpha_y(b)(yv)_+ \bigr]\\
 \nonumber  &= (xs)_-\phi\alpha_x(a)\phi\alpha_x(p_t)\phi\alpha_y(b)(yv)_+ \\
 \nonumber  &=
 (xs)_-\phi\bigl(\alpha_x(ap_t)\alpha_y(b)\bigr)(yv)_+. \qquad\qquad\Box
\end{align}

\begin{corollary}\label{sum} $R= \sum_{s\in S,\, t\in T} s_-\phi(A)t_+=
\sum_{s\in S,\, t\in T} s_-\phi(p_sAp_t)t_+$. \end{corollary}

\begin{proof} The second equality is clear from Lemma
\ref{relations1}(d). Let $R'$ denote the sum in question. Clearly
$R'$ is closed under addition, and it is closed under
multiplication by Lemma \ref{relations1}(f). Also, $1_R=
1_-\phi(1_A)1_+ \in R'$. Thus, $R'$ is a unital subring of $R$.

Since the images of $\phi$, $s\mapsto s_-$, and $t\mapsto t_+$ are
contained in $R'$, we can view these as maps into $R'$. The
universal property for $R$ then implies that there is a unique
unital ring homomorphism $\psi: R\rightarrow R'$ such that
$\psi\phi= \phi$ while $\psi(s_-)= s_-$ for $s\in S$ and
$\psi(t_+)= t_+$ for $t\in T$. Consequently, $\psi$ acts as the
identity on $R'$, whence $\psi(R)= R'$. Moreover, if we view
$\psi$ as a ring homomorphism $R\rightarrow R$, we have $\psi\phi=
\id_R\phi$ while $\psi(s_-)= \id_R(s_-)$ for $s\in S$ and
$\psi(t_+)= \id_R(t_+)$ for $t\in T$. Now the universal property
for $R$ implies that $\psi= \id_R$, and therefore $R= \psi(R)=
R'$. \end{proof}

\begin{lemma}\label{idealfact} {\rm (a)} The set $I= \bigcup_{s\in S}
\ker(\alpha_s)$ is an ideal of $A$, contained in $\ker(\phi)$.

{\rm (b)} $\alpha_s^{-1}(I)= I$ for all $s\in S$, and $\alpha_t(I)
\subseteq I$ for all $t\in T$.

{\rm (c)} $\alpha$ induces a monoid homomorphism $\alpha':
T\rightarrow \ndz(A/I)$, and $\alpha'_s$ is injective for all
$s\in S$.

{\rm (d)} $\sop*_\alpha A*_\alpha T= \sop*_{\alpha'}
(A/I)*_{\alpha'} T$. \end{lemma}

\begin{proof} (a) If $s_1,s_2\in S$, there exist $s\in S$ and $t\in
T$ such that $ss_1= ts_2$. Then $ss_1\in S$ and
$\ker(\alpha_{s_i}) \subseteq \ker(\alpha_{ss_1})$ for $i=1,2$.
Thus, the ideals $\ker(\alpha_s)$ for $s\in S$ are upward directed
under inclusion, whence $I$ is an ideal of $A$. Each
$\ker(\alpha_s) \subseteq \ker(\phi)$ by Lemma
\ref{relations1}(b), and so $I\subseteq \ker(\phi)$.

(b) If $t\in T$ and $s\in S$, there exist $s'\in S$ and $t'\in T$
such that $s't= t's$. Then $\alpha_{s'}\alpha_t(\ker(\alpha_s))
=0$, and so $\alpha_t(\ker(\alpha_s)) \subseteq \ker(\alpha_{s'})
\subseteq I$. This shows that $\alpha_t(I) \subseteq I$ for all
$t\in T$.

Now if $s\in S$, the previous paragraph implies that $I \subseteq
\alpha_s^{-1}(I)$. If $a\in \alpha_s^{-1}(I)$, then
$\alpha_s(a)\in \ker(\alpha_{s'})$ for some $s'\in S$, whence
$a\in \ker(\alpha_{s's}) \subseteq I$. Therefore $\alpha_s^{-1}(I)
=I$.

(c)(d) These are clear from (a) and (b). \end{proof}

\begin{proposition}
\label{gradedring} The ring $R$ has an $\sinvt$-grading $R=
\bigoplus_{x\in\sinvt} R_x$ where each $R_x= \bigcup_{s^{-1}t=x}
s_-\phi(A)t_+$. \end{proposition}

\begin{proof} We can view $R$ as a left $A$-module via $\phi$, and
the relations in $R$ imply that each $s_-\phi(A)t_+$ is a left
$A$-submodule. If $s_1,s_2\in S$ and $t_1,t_2\in T$ such that
$s_1^{-1}t_1= s_2^{-1}t_2$, there exist $u_1,u_2\in S$ such that
$u_1s_1= u_2s_2$ and $u_1t_1= u_2t_2$, whence Lemma
\ref{relations1}(e) implies that $(s_i)_-\phi(A)(t_i)_+ \subseteq
(u_1s_1)_-\phi(A)(u_1t_1)_+$ for $i=1,2$. Thus, each $R_x$ is a
directed union of left $A$-submodules of $R$, and so is a left
$A$-submodule itself.

It is clear from Corollary \ref{sum} that $R= \sum_{x\in\sinvt}
R_x$, and from Lemma \ref{relations1}(f) that $R_xR_y \subseteq
R_{xy}$ for all $x,y\in\sinvt$. Hence, it only remains to show
that the sum of the $R_x$ is a direct sum.

Let $R'$ denote the external direct sum of the $R_x$, and set $E'=
\ndz(R')$. There is a unital ring homomorphism $\lambda:
A\rightarrow E'$ such that each $\lambda(a)$ is the left
$A$-module multiplication by $a\in A$.

Given $s\in S$, observe that $s_-R_x \subseteq R_{s^{-1}x}$ for
all $x\in\sinvt$. Hence, there exists $\mu_s\in E'$ such that
$\mu_s(b)_y= s_-b_{sy}$ for all $b\in R'$ and $y\in\sinvt$. Since
$\phi(a)s_-= s_-\phi\alpha_s(a)$ for $a\in A$, we see that
$\lambda(a)\mu_s= \mu_s\lambda\alpha_s(a)$ for $a\in A$. Observe
also that $s\mapsto \mu_s$ is a monoid homomorphism
$\sop\rightarrow E'$.

Given $t\in T$, it follows from Lemma \ref{relations1}(f) that
$t_+R_x \subseteq R_{tx}$ for all $x\in\sinvt$. Hence, there
exists $\nu_t\in E'$ such that $\nu_t(b)_y= \sum_{tx=y} t_+b_x$
for $b\in R'$ and $y\in\sinvt$. Since $t_+\phi(a)=
\phi\alpha_t(a)t_+$ for $a\in A$, we see that $\nu_t\lambda(a)=
\lambda\alpha_t(a)\nu_t$ for $a\in A$. Observe also that $t\mapsto
\nu_t$ is a monoid homomorphism $T\rightarrow E'$.

Since $s_-s_+=1$ and $s_+s_-= \phi(p_s)$ for $s\in S$, we see that
$\mu_s\nu_s= \id_{R'}= 1_{E'}$ and $\nu_s\mu_s= \lambda(p_s)$ for
$s\in S$. Now by the universal property of $R$, there exists a
unital ring homomorphism $\psi: R\rightarrow E'$ such that
$\psi\phi= \lambda$ while $\psi(s_-)= \mu_s$ for $s\in S$ and
$\psi(t_+)= \nu_t$ for $t\in T$.

Note that $1_R= 1_-\phi(1)1_+ \in R_1$, so there exists $e\in R'$
such that $e_1=1$ while $e_z=0$ for all $z\ne 1$. Given $s\in S$,
$a\in A$, and $t\in T$, we observe that
$$\bigl[ \psi(s_-\phi(a)t_+)(e) \bigr]_{s^{-1}t}= \bigl[
\mu_s\lambda(a)\nu_t(e) \bigr]_{s^{-1}t}= s_-\phi(a)t_+$$
 and all
other components of $\psi(s_-\phi(a)t_+)(e)$ are zero. Hence, for
$x\in\sinvt$ and $b\in R_x$, we have $[\psi(b)(e)]_x=b$ while
$[\psi(b)(e)]_y=0$ for all $y\ne x$. Consequently, if $b_1+\cdots
+b_n=0$ for some $b_i\in R_{x_i}$ where the $x_i$ are distinct
elements of $\sinvt$, then $b_i=
[\psi(b_1+\cdots+b_n)(e)]_{x_i}=0$ for all $i$. Therefore
$\sum_{x\in\sinvt} R_x= \bigoplus_{x\in\sinvt} R_x$, as desired.
\end{proof}

To completely pin down the elements of $R$, we need to know the
relations holding in each homogeneous component $R_x$. In
particular, if $p_sap_t\in \ker(\phi)$, then $s_-\phi(a)t_+ =0$ by
Lemma \ref{relations1}(d), and we would like to show that
$s_-\phi(a)t_+ =0$ only when $p_sap_t\in \ker(\phi)$. For this
purpose, we set up another representation of $R$ on a left
$A$-module.

\begin{lemma}\label{moduleact} Let $u,s\in S$ and $t\in T$.

{\rm (a)} The map $*: A\times p_sAp_t \rightarrow p_sAp_t$ given
by the rule $a*b := \alpha_s(a)b$ turns the abelian group
$p_sAp_t$ into a left $A$-module.

{\rm (b)} The restriction of $\alpha_u$ to $p_sAp_t$ is a left
$A$-module homomorphism $p_sAp_t \rightarrow p_{us}Ap_{ut}$.
\end{lemma}

\begin{proof} Part (a) is clear because $\alpha_s$ is a unital
ring homomorphism from $A$ to $p_sAp_s$, while part (b) follows
because $\alpha_{us}= \alpha_u\alpha_s$.
\end{proof}

Each homogeneous component $R_x$ of $R$ turns out to be a direct
limit of the rectangular corners $p_sAp_t$ over pairs $(s,t)$ such
that $s^{-1}t=x$. However, there is no natural partial order on
the set of these pairs -- the limit has to be taken over a small
category.

\begin{definition}\label{directsyst} {\rm For $x\in\sinvt$, let
$\D_x$ be the small category in which the objects are all pairs
$(s,t)\in S\times T$ such that $s^{-1}t=x$, the morphisms from an
object $(s,t)$ to an object $(s',t')$ are those elements $u\in S$
such that $us=s'$ and $ut=t'$, and composition of morphisms is
given by the multiplication in $S$. The Ore and saturation
conditions on $S$ imply that $\D_x$ is directed: given any objects
$(s_1,t_1)$ and $(s_2,t_2)$ in $\D_x$, there exist an object
$(s,t)$ and morphisms $u_i: (s_i,t_i) \rightarrow (s,t)$ in $\D_x$
for $i=1,2$. Consequently, colimits based on $\D_x$ are directed
colimits.

Taking account of Lemma \ref{moduleact}, there is a functor $F_x:
\D_x \rightarrow \amod$ such that $F_x(s,t)= p_sAp_t$ for all
objects $(s,t)$ in $\D_x$ and $F_x(u)= \alpha_u|_{p_sAp_t}$ for
all morphisms $u: (s,t) \rightarrow (us,ut)$ in $\D_x$. Let $M_x$
denote the colimit of $F_x$, with natural maps $\eta_{s,t}:
p_sAp_t \rightarrow M_x$ for objects $(s,t)$ in $\D_x$. Since
$M_x$ is a directed colimit, it is the union of its submodules
$\eta_{s,t}(p_sAp_t)$ for $(s,t) \in \D_x$. Note that if $b_i \in
p_{s_i}Ap_{t_i}$ for $i=1,2$, where $(s_i,t_i) \in D_x$, then
$\eta_{s_1,t_1}(b_1)= \eta_{s_2,t_2}(b_2)$ if and only if there
exist $u_1,u_2 \in S$ such that $u_1s_1= u_2s_2$ and $u_1t_1=
u_2t_2$ while also $\alpha_{u_1}(b_1)= \alpha_{u_2}(b_2)$.}
\end{definition}

\begin{lemma}\label{trickymap} Let $s\in S$, $t\in T$, and $x\in\sinvt$.

{\rm (a)} There exists an additive map $\sigma_s: M_x\rightarrow
M_{s^{-1}x}$ such that $\sigma_s\eta_{u,v}(b)=
\eta_{us,v}(p_{us}b)$ for $u^{-1}v=x$ and $b\in p_uAp_v$.

{\rm (b)} $a\sigma_s(m)= \sigma_s(\alpha_s(a)m)$ for $a\in A$ and
$m\in M_x$.

{\rm (c)} There exists an additive map $\tau_t: M_x\rightarrow
M_{tx}$ such that $\tau_t\eta_{u,v}(b)= \eta_{w,zv}\alpha_z(b)$
for $u^{-1}v=x$, $b\in p_uAp_v$, and $w\in S$, $z\in T$ such that
$wt=zu$.

{\rm (d)} $\tau_t(am)= \alpha_t(a)\tau_t(m)$ for $a\in A$ and
$m\in M_x$.
\end{lemma}

\noindent\emph{Proof}. (a) For each $(u,v)\in \D_x$, we have
$(us,v)\in \D_{s^{-1}x}$, and there is an additive map $p_uAp_v
\rightarrow M_{s^{-1}x}$ given by $b\mapsto \eta_{us,v}(p_{us}b)$.
Moreover, if $w\in S$ then $\eta_{wus,wv}(p_{wus}\alpha_w(b))=
\eta_{wus,wv} \alpha_w(p_{us}b)= \eta_{us,v}(p_{us}b)$. Thus, our
maps to $M_{s^{-1}x}$ are compatible with the functor $F_x$, and
so there exists a unique additive map $\sigma_s$ as described.

(b) If $m=\eta_{u,v}(b)$ for $u$, $v$, $b$ as in (a), then
\begin{align}
\nonumber  a\sigma_s(m) &= a\eta_{us,v}(p_{us}b) =
\eta_{us,v}(a*(p_{us}b)) = \eta_{us,v}(\alpha_{us}(a)p_{us}b) =
\eta_{us,v}(p_{us}\alpha_{us}(a)b) \\
\nonumber  &=\eta_{us,v}(p_{us}(\alpha_s(a)*b)) =
\sigma_s\eta_{u,v}(\alpha_s(a)*b) =\sigma_s(\alpha_s(a)m).
 \end{align}

 (c) Fix $(u,v)\in \D_x$, choose $w\in S$, $z\in T$ such that
 $wt=zu$, and note that $tx=w^{-1}zv$. Since $\alpha_z(p_uAp_v)
 \subseteq p_{zu}Ap_{zv} \subseteq p_wAp_{zv}$, the composition of
 $\eta_{w,zv}$ with the restriction of $\alpha_z$ to $p_uAp_v$
 gives an additive map $p_uAp_v \rightarrow M_{tx}$. Suppose also
 $w_1\in S$ and $z_1\in T$ such that $w_1t=z_1u$. Then
 $w_1^{-1}z_1= tu^{-1}= w^{-1}z$, so there exist $r_1,r\in S$ such
 that $r_1w_1=rw$ and $r_1z_1= rz$. Since also $r_1z_1v= rzv$ and
 $\alpha_{r_1}\alpha_{z_1}= \alpha_r\alpha_z$, it follows that
 $\eta_{w_1,z_1v}\alpha_{z_1}= \eta_{w,zv}\alpha_z$ on $p_uAp_v$.
 Thus, we obtain a well-defined additive map $f_{u,v}: p_uAp_v
 \rightarrow M_{tx}$ which agrees with $\eta_{w,zv}\alpha_z$ for
 any $w\in S$ and $z\in T$ with $wt=zu$.

 Now consider a morphism $r: (u,v) \rightarrow (ru,rv)$ in $\D_x$.
 There exist $w\in S$ and $z\in T$ such that $wt= z(ru)$, so that
 $f_{ru,rv}$ is given by $\eta_{w,zrv}\alpha_z$. Since $wt=
 (zr)u$, we also have that $f_{u,v}$ is given by
 $\eta_{w,zrv}\alpha_{zr}$, and so $f_{u,v}$ equals the
 composition of $f_{ru,rv}$ with the restriction of $\alpha_r$ to
 $p_uAp_v$. Thus, the maps $f_{.,.}$ are compatible with $F_x$,
 and so there exists a unique additive map $\tau_t$ as described.

 (d) If $m= \eta_{u,v}(b)$ with $u$, $v$, $b$, $w$, $z$ as in (c),
 then
\begin{align}
\nonumber \tau_t(am) &= \tau_t\eta_{u,v}(a*b)=
 \tau_t\eta_{u,v}(\alpha_u(a)b)=
 \eta_{w,zv}\alpha_z(\alpha_u(a)b)=
  \eta_{w,zv}(\alpha_w\alpha_t(a)\alpha_z(b)) \\
\nonumber  &= \eta_{w,zv}(\alpha_t(a)*\alpha_z(b))=
  \alpha_t(a)\eta_{w,zv}\alpha_z(b)= \alpha_t(a)\tau_t(m).
  \qquad\qquad\square
  \end{align}

\begin{proposition}\label{isomodule} For each $x\in \sinvt$, there
is a left
$A$-module isomorphism $\theta_x: M_x\rightarrow R_x$ such that
$\theta_x\eta_{u,v}(b)= u_-\phi(b)v_+$ for $u^{-1}v=x$ and $b\in
p_uAp_v$.
\end{proposition}

\begin{proof} In view of Lemma \ref{relations1}(e), for each
$x\in\sinvt$ there is
a unique additive map $\theta_x: M_x\rightarrow R_x$ as described.
If $m= \eta_{u,v}(b)$ with $u$, $v$, $b$ as above, then for $a\in
A$ we have
$$\theta_x(am) = \theta_x\eta_{u,v}(a*b)=
\theta_x\eta_{u,v}(\alpha_u(a)b)
 = u_-\phi\alpha_u(a)\phi(b)v_+=
\phi(a)u_-\phi(b)v_+=
 a\theta_x(m).$$
 Thus, $\theta_x$ is a left $A$-module homomorphism.
It is surjective by definition of $R_x$, and so it only remains to
show that $\ker(\theta_x) =0$.

Form the left $A$-module $M := \bigoplus_{x\in \sinvt} M_x$, set
$E= \ndz(M)$, and for each $a\in A$ let $\lambda(a) \in E$ be the
map given by left multiplication by $a$. Then we have a unital
ring homomorphism $\lambda: A\rightarrow E$.

For all $x\in\sinvt$, use the same notations $\sigma_s$ and
$\tau_t$ for the additive maps $M_x\rightarrow M_{s^{-1}x}$ and
$M_x\rightarrow M_{tx}$ described in Lemma \ref{trickymap}, and
also for the corresponding homogeneous maps on $M$. Thus, for
$s\in S$ and $t\in T$ we have additive maps $\sigma_s,\tau_t\in E$
such that $\sigma_s(m)_y= \sigma_s(m_{sy})$ and $\tau_t(m)_y=
\sum_{tx=y} \tau_t(m_x)$ for $m\in M$ and $y\in\sinvt$. Lemma
\ref{trickymap} also shows that $\lambda(a)\sigma_s=
\sigma_s\lambda\alpha_s(a)$ and $\tau_t\lambda(a)=
\lambda\alpha_t(a)\tau_t$ for $a\in A$.

It is easily checked that $s\mapsto \sigma_s$ and $t\mapsto
\tau_t$ are monoid homomorphisms $\sop\rightarrow E$ and
$T\rightarrow E$. Now consider $m=\eta_{u,v}(b)\in M_x$ for $x$,
$u$, $v$, $b$ as in Lemma \ref{trickymap}. There exist $w\in S$
and $z\in T$ such that $ws=zu$, and
\begin{align}
\nonumber \sigma_s\tau_s(m) &= \sigma_s\eta_{w,zv}\alpha_z(b)=
\eta_{ws,zv}(p_{ws}\alpha_z(b))= \eta_{zu,zv}(p_{zu}\alpha_z(b))
\\
\nonumber  &= \eta_{zu,zv}\alpha_z(p_ub)= \eta_{u,v}(b)= m.
 \end{align}
 It follows that $\sigma_s\tau_s= 1_E$ in $E$. Next, note that
 $u\in S$ and $1\in T$ with $u\cdot s= 1\cdot us$. Hence,
 $$\tau_s\sigma_s(m)= \tau_s\eta_{us,v}(p_{us}b)=
 \eta_{u,v}\alpha_1(p_{us}b)= \eta_{u,v}(p_s*b)= p_sm.$$
 It follows that $\tau_s\sigma_s=\lambda(p_s)$ in $E$.

 By the universal property of $R$, there is a unital ring
 homomorphism $\psi: R\rightarrow E$ such that $\psi\phi= \lambda$
 while $\psi(s_-)= \sigma_s$ for $s\in S$ and $\psi(t_+)=
 \tau_t$ for $t\in T$.

 Define $e\in M$ so that $e_1= \eta_{1,1}(1)$ while $e_z=0$ for
 all $z\ne 1$. We claim that $[(\psi\theta_x(m))(e)]_x =m$ for
 $x\in\sinvt$ and $m\in M_x$. Write $m= \eta_{u,v}(b)$ where
 $u^{-1}v=x$ and $b\in p_uAp_v$. Then $\psi\theta_x(m)=
 \psi(u_-\phi(b)v_+)= \sigma_u\lambda(b)\tau_v$ and so
 \begin{align}
 \nonumber [(\psi\theta_x(m))(e)]_x &= \sigma_u\lambda(b)\tau_v\eta_{1,1}(1)=
 \sigma_u\lambda(b)\eta_{1,v}\alpha_v(1)=
 \sigma_u\eta_{1,v}(b*p_v) \\
 \nonumber  &= \sigma_u\eta_{1,v}(b)= \eta_{u,v}(p_ub)= \eta_{u,v}(b)= m,
  \end{align}
  as claimed.

  The claim immediately implies that $\ker(\theta_x)=0$ for all
  $x\in\sinvt$, as desired. \end{proof}

As the reader will note, the grading $R= \bigoplus_{x\in\sinvt}
R_x$ can also be obtained from the proof of Proposition
\ref{isomodule}, and so Proposition \ref{gradedring} could have
been omitted. However, we think the latter proposition is helpful
in orienting the reader.

 \begin{corollary}\label{kernel} Let $s\in S$, $t\in T$, and $a\in A$.
  Then $s_-\phi(a)t_+=0$ if and only if $p_sap_t\in
  \ker(\alpha_{s'})$ for some $s'\in S$. In particular,
 $\ker(\phi)= \bigcup_{s'\in S} \ker(\alpha_{s'})$. \end{corollary}

 \begin{proof} By Lemma \ref{relations1}(d), $s_-\phi(a)t_+= s_-\phi(b)t_+$ where
 $b= p_sap_t$. Then Proposition \ref{isomodule} yields $\theta_x\eta_{s,t}(b)=
 s_-\phi(a)t_+$ where $x= s^{-1}t$. Since $\theta_x$ is an
 isomorphism, $s_-\phi(a)t_+=0$ if and only if $\eta_{s,t}(b)=0$,
 which happens if and only if $\alpha_{s'}(b)=0$ for some $s'\in
 S$. This verifies the first statement of the corollary. The
 second follows on taking $s=t=1$. \end{proof}

\begin{noname}\label{alphasinject} {\rm
 As Lemma \ref{idealfact} shows, we can always reduce to the case where
 $\alpha_s$ is injective for all $s\in S$. In that case, $\phi$ is
 injective by Corollary \ref{kernel}, and so we can identify $A$ with the
 unital subring $\phi(A) \subseteq R$. All of the relations in $R$
 simplify in this case:

\begin{enumerate}

\item $t_+a= \alpha_t(a)t_+$ for all $a\in A$ and $t\in T$;

\item $as_-= s_-\alpha_s(a)$ for all $a\in A$ and $s\in S$;

\item $s_-s_+= 1$ for all $s\in S$;

\item $s_+s_-= p_s$ for all $s\in S$;

\item $R$ has an $\sinvt$-grading $R= \bigoplus_{x\in\sinvt} R_x$
where each $R_x= \bigcup_{s^{-1}t=x} s_-At_+$;

\item $s_-at_+= s_-p_sap_tt_+$ for $s\in S$, $t\in T$, and $a\in
A$, and $s_-at_+=0$ if and only if $p_sap_t=0$;

\item Let $x=s_1^{-1}t_1= s_2^{-1}t_2\in\sinvt$ for some
$s_1,s_2\in S$, $t_1,t_2\in T$, and let $a_1,a_2\in A$. Then
$(s_1)_-a_1(t_1)_+= (s_2)_-a_2(t_2)_+$ if and only if there exist
$u_1,u_2\in S$ such that $u_1s_1=u_2s_2$ and $u_1t_1=u_2t_2$ while
also $\alpha_{u_1}(p_{s_1}a_1p_{t_1})=
\alpha_{u_2}(p_{s_2}a_2p_{t_2})$.
\end{enumerate} }
\end{noname}

\section{The case $S=T=\Zplus$; Examples}

\begin{noname}\label{data2} {\rm
 For the remainder of the paper,
we take advantage of Lemma \ref{idealfact} and assume that
$\alpha_s$ is injective for all $s\in S$. Thus, the relations in
$R=\sopat$ take the simplified form given in (\ref{alphasinject}).
Moreover, we assume that the maps $\alpha_s$ are \emph{corner
isomorphisms}, that is, each $\alpha_s$ is an isomorphism of $A$
onto $p_sAp_s$. Finally, we assume that $S=T$ is a submonoid of a
group $G$ which is its group of left fractions, that is,
$G=\sinvs$. These conventions are to remain in effect for the rest
of the paper. }
\end{noname}

\begin{noname}\label{Gtotord} {\rm
A particularly nice setting is the case when $G$ is a left totally
ordered group with positive cone $G^+=S$ (thus $G=S^{-1}\cup S$
and $S^{-1}\cap S= \{1\}$). In this case, the elements of $R$ can
be expressed in a simpler way, namely in the form $\sum _{s\in
S}s_{-}a_s+ \sum_{t\in S}a_tt_+$. To achieve this, we need to be
able to simplify individual terms $s_-at_+$, for $s,t\in S$ and
$a\in A$. If $s\le t$, then $s^{-1}t\ge 1$, whence $u:= s^{-1}t\in
S$. Then $s_-at_+= s_-a(su)_+= s_-p_sap_ss_+u_+$. Because of our
current convention that $\alpha_s: A\rightarrow p_sAp_s$ is an
isomorphism, $p_sap_s=\alpha_s(b)$ for some $b\in A$, and
therefore $s_-at_+= s_-\alpha_s(b)s_+u_+= bs_-s_+u_+= bu_+$. On
the other hand, if $s\ge t$, then $v:= t^{-1}s\in S$ and $s_-at_+=
v_-c$ where $c=\alpha_t^{-1}(p_tap_t)$. }
\end{noname}

\begin{noname}\label{GequalZ} {\rm
We now specialize to the case where $S$ is the additive monoid
$\Zplus$, so that $G=\Z$. Here the monoid homomorphism $\alpha:
S\rightarrow \ndr(A)$ is determined by $\alpha_1$, and so we
change notation, writing $\alpha$ and $p$ for $\alpha_1$ and
$p_1$. Thus, $\alpha$ is now an isomorphism $A\rightarrow pAp$,
and the monoid homomorphism $S\rightarrow \ndr(A)$ is given by the
rule $n\mapsto \alpha^n$. Let $t$ denote the generator $1\in
\Zplus=S$. Since the maps $s\mapsto s_\pm$ are monoid
homomorphisms into the multiplicative structure of $R$, we have
$n_\pm= (t_\pm)^n =: t^n_\pm$ for $n\in \Zplus$, and
$$at^n_-= t^n_-\alpha^n(a) \qquad\quad \text{and} \qquad\quad
t^n_+a= \alpha^n(a)t^n_+$$
 for all $a\in A$ and $n\in \Zplus$.

In view of (\ref{Gtotord}), the elements $r\in R= \zpazp$ can all
be written as `polynomials' of the form
$$r=a_{n}t_+^n+\ldots
+a_{1}t_++a_0+t_{-}a_{-1}+\ldots t_{-}^m a_{-m},$$
 with coefficients $a_i\in A$. Because of this similarity of $R$ with
a skew-Laurent polynomial ring, we shall use the notation $R=
\atpmalpha$. Proposition \ref{gradedring} shows
 that $R$ is a $\Z$-graded ring $R= \bigoplus_{i\in\Z} R_i$, and
from the discussion above we see that $R_i= At^i_+$ for $i>0$ and
$R_i= t_-^{-i}A$ for $i<0$, while $A_0=A$.

Our construction of $\zpazp$ is an exact algebraic analog of the
construction of the crossed product of a C*-algebra by an
endomorphism introduced by Paschke \cite{Paschke}. In fact, if $A$
is a C*-algebra and the corner isomorphism $\alpha$ is a
*-homomorphism, then Paschke's C*-crossed product, which he
denotes $A\rtimes _{\alpha}\N$, is just the completion of $\zpazp$
in a suitable norm. }
\end{noname}

Note again that any ring $R= \atpmalpha$ is $\Z$-graded, with $A=
R_0$. Moreover, $t_+$ is a left invertible element of $R_1$ with a
particular left inverse $t_-\in R_{-1}$, and $\alpha$ can be
recovered from the rule $\alpha(a)= t_+at_-$. These observations
allow us to recognize rings of the form $\atpmalpha$ among
$\Z$-graded rings, as follows.

\begin{lemma}\label{recognize}
Let $D= \bigoplus_{i\in\Z} D_i$ be a $\Z$-graded ring containing
elements $t_+\in D_1$ and $t_-\in D_{-1}$ such that $t_-t_+=1$.
Then there is a corner isomorphism $\alpha: D_0\rightarrow
t_+t_-D_0t_+t_-$ given by the rule $\alpha(d)= t_+dt_-$, and $D=
\dztpmalpha$.
\end{lemma}

\begin{proof} It is clear that $t_+t_-$ is an idempotent in $D_0$,
and that the given rule defines an isomorphism $\alpha:
D_0\rightarrow t_+t_-D_0t_+t_-$. Hence, there exists a fractional
skew monoid ring $\Dtil= D_0[\ttil_+,\ttil_-;\alpha]$. Since
$t_+d= \alpha(d)t_+$ and $dt_-= t_-\alpha(d)$ for all $d\in D$,
the identity map on $D_0$ extends uniquely to a ring homomorphism
$\phi: \Dtil \rightarrow D$ such that $\phi(\ttil_\pm)= t_\pm$. It
remains to show that $\phi$ is an isomorphism. Note that since
$t_+^i \in D_i$ and $t_-^i \in D_{-i}$ for all $i\in\N$, the map
$\phi$ is a homomorphism of graded rings. Thus, we need only show
that $\phi$ maps each homogeneous component $\Dtil_i$
isomorphically onto $D_i$. This is already given when $i=0$.

Now let $i>0$. If $x\in \Dtil_i$, then $x= d\ttil_+^i$ for some
$d\in D_0$, and $\phi(x)= dt_+^i$. If $\phi(x)=0$, then
$d\alpha^i(1)= dt_+^it_-^i= 0$ in $D_0$, whence $x=
d\alpha^i(1)\ttil_+^i =0$ in $\Dtil$. Thus, the restriction of
$\phi$ to $\Dtil_i$ is injective. Further, if $y\in D_i$, then
$yt_-^i\in D_0$ and $\phi \bigl( (yt_-^i)\ttil_+^i \bigr)=
yt_-^it_+^i= y$. Therefore $\phi$ maps $\Dtil_i$ isomorphically
onto $D_i$. A symmetric argument shows that this also holds for
$i<0$, completing the proof. \end{proof}

\begin{example}\label{algCK}
 {\rm \emph{An algebraic version of the Cuntz-Krieger
algebras}. We give an algebraic version of the C*-algebras
$\mathcal{O}_A$ introduced in \cite{c-k} (now called
``Cuntz-Krieger algebras'' in the literature), and show that they
may be expressed in the form $\btpmalpha$ for ultramatricial
algebras $B$ and proper corner isomorphisms $\alpha$. The latter
statement is parallel to the corresponding C*-algebra result:
$\mathcal{O}_A= B\rtimes_\alpha \N$ for a suitable approximately
finite dimensional C*-algebra $B$ (essentially in \cite{c-k};
discussed explicitly in \cite[Example 2.5]{rordamclassif}).

Let $k$ be an arbitrary field and $A= (a_{ij})$ an $n\times n$
matrix over $k$, with $a_{ij}\in \{0,1\}$ for all $i,j$. To avoid
degenerate and trivial cases, we assume that no row or column of
$A$ is identically zero, and that $A$ is not a permutation matrix.
We define the \emph{algebraic Cuntz-Krieger algebra associated to
$A$} to be the $k$-algebra $C=\mathcal{CK}_A(k)$ with generators
$x_1,y_1,\dots,x_n,y_n$ and relations
 \begin{enumerate}
\item $x_iy_ix_i=x_i$ and $y_ix_iy_i=y_i$ for all $i$;

\item $x_iy_j=0$ for all $i\ne j$;

\item $x_iy_i= \sum_{j=1}^n a_{ij}y_jx_j$ for all $i$;

\item $\sum_{j=1}^n y_jx_j=1$.
 \end{enumerate}
Note that all the $x_iy_i$ and $y_jx_j$ are idempotents, and that
the $y_jx_j$ are pairwise orthogonal. The free algebra $k\langle
X_1,Y_1, \dots, X_n,Y_n\rangle$ can be given a $\Z$-grading in
which the $X_i$ have degree $-1$ while the $Y_i$ have degree $1$,
and the relators $X_iY_iX_i-X_i$ etc. corresponding to (1)--(4)
are all homogeneous. Hence, $C$ inherits a $\Z$-grading such that
each $x_i\in C_{-1}$ and each $y_i\in C_1$.

Now set $N=\{1,\dots,n\}$. Given $\mu= (\mu_1,\dots,\mu_\ell)\in
N^\ell$ for some $\ell$, we set $x_\mu= x_{\mu_1} x_{\mu_2} \cdots
x_{\mu_\ell}$ and $y_\mu= y_{\mu_1} y_{\mu_2} \cdots
y_{\mu_\ell}$. The case $\ell=0$ is allowed, with the conventions
that $N^0=\{\varnothing\}$ and $x_\varnothing= y_\varnothing= 1$.
The subalgebra $B=C_0$ of $C$ is the $k$-linear span of the set
 $$\{y_\mu x_\nu \mid \mu,\nu\in N^\ell,\ \ell\in\Zplus\}.$$
As in \cite[Proposition 2.3 and following discussion]{c-k}, $B$ is
an ultramatricial $k$-algebra, and $K_0(B)$ is isomorphic (as an
ordered group) to the direct limit of the sequence
 $$\Z^n \overset{A}\longrightarrow \Z^n \overset{A}\longrightarrow
 \Z^n \overset{A}\longrightarrow \cdots,$$
with the class $[B]\in K_0(B)$ corresponding to the image of the
order-unit $(1,1,\dots,1)^{\text{tr}}$ in the first $\Z^n$. (See
\cite[Chapter 15]{vnrr} for a development of ultramatricial
algebras and their classification via $K_0$.)

For $i=1,\dots,n$, let $e_i$ denote the sum of those $y_jx_j$ for
which $y_jx_j \le x_iy_i$ but $y_jx_j \not\le x_my_m$ for any
$m<i$. These $e_i$ are pairwise orthogonal idempotents in $B$,
with each $e_i\le x_iy_i$. Since the matrix $A$ has no identically
zero columns, each $y_jx_j$ lies below some $x_iy_i$, and so each
$y_jx_j$ lies below some $e_i$. In fact, $y_jx_j\le e_i$ where $i$
is the least index such that $a_{ij}=1$. From relation (4), it
follows that $\sum_{i=1}^n e_i=1$. Next, note that the elements
$y_ie_ix_i$ are pairwise orthogonal idempotents in $B$ (because
$e_ix_iy_i= e_i$ for all $i$), whence the sum $p:= y_1e_1x_1+
\cdots+ y_ne_nx_n$ is an idempotent in $B$. Moreover,
$x_ip=e_ix_i$ and $py_i= y_ie_i$ for all $i$. We claim that $p\ne
1$.

If $p=1$, then each $x_i=e_ix_i$, whence each $x_iy_i=e_i$. Then
the $x_iy_i$ are pairwise orthogonal. In view of the relations
(3), it follows that each column of $A$ has only one nonzero
entry. Since $A$ has no identically zero rows, it must be a
permutation matrix, contradicting our assumptions. Therefore $p\ne
1$, as claimed.

Now set $t_-= e_1x_1+ \cdots+ e_nx_n\in C_{-1}$ and $t_+= y_1e_1+
\cdots+ y_ne_n\in C_1$. Then $t_+t_-=p$, and
 $$t_-t_+= \sum_{i=1}^n e_ix_iy_ie_i= \sum_{i,j=1}^n
 a_{ij}e_iy_jx_je_i= \sum_{j=1}^n y_jx_j =1,$$
because each $y_jx_j\le e_i$ for precisely one $i$, and $a_{ij}=1$
for that $i$. Hence, there is a proper corner isomorphism $\alpha:
B\rightarrow pBp$ given by the rule $\alpha(b)= t_+bt_-$, and we
conclude from Lemma \ref{recognize} that
 $$C=\mathcal{CK}_A(k)= \btpmalpha. \qquad\qquad\square$$  }
\end{example}

In case the matrix $A$ in Example \ref{algCK} has all of its
entries equal to $1$, the relations for the algebra
$\mathcal{CK}_A(k)$ reduce to
 \begin{enumerate}
\item $x_iy_j= \delta_{i,j}$ for all $i,j$;

\item $\sum_{j=1}^n y_jx_j=1$.
 \end{enumerate}
Thus in this case, $\mathcal{CK}_A(k)$ is the \emph{Leavitt
algebra} $V_{1,n}(k)$ first studied in \cite{Lvtt}. (The notation
$V_{1,n}$ was introduced in \cite{Berg}.) There is a related
Leavitt algebra $U_{1,n}(k)$ which, as we now show, can also be
presented as a fractional skew monoid ring.

\begin{example}\label{U1n} {\rm
 Let $k$ be a field and $n\in\N$. The algebra $U=U_{1,n}(k)$ is
 the $k$-algebra with generators
$x_1,y_1,\dots,x_n,y_n$ and relations $x_iy_j=\delta_{i,j}$ for
all $i,j$. (Thus, $V_{1,n}(k)$ is the factor algebra of
$U_{1,n}(k)$ modulo the ideal generated by $1-\sum_{j=1}^n
y_jx_j$.) The elements $y_1x_1,\dots,y_nx_n$ are pairwise
orthogonal idempotents in $U$. As in Example \ref{algCK}, there is
a $\Z$-grading on $U$ such that each $x_i\in U_{-1}$ and each
$y_i\in U_1$.

Set $N=\{1,\dots,n\}$ and define $x_\mu,y_\mu\in U$ for $\mu\in
N^\ell$ as in Example \ref{algCK}. In $U$, the set
 $$\{ y_\mu x_\nu \mid \mu\in N^\ell,\ \nu\in N^m,\ \ell,m\in\Zplus
 \}$$
forms a $k$-basis. We again set $B=U_0$, which is the $k$-linear
span of the set
 $$\{y_\mu x_\nu \mid \mu,\nu\in N^\ell,\ \ell\in\Zplus\},$$
and as before, $B$ is ultramatricial. It is isomorphic to a direct
limit of the algebras
 $$M_{n^i}(k)\times M_{n^{i-1}}(k)\times \cdots\times M_n(k)\times k,$$
the ordered group $K_0(B)$ is isomorphic to the direct limit of a
sequence $\Z \rightarrow \Z^2 \rightarrow \Z^3 \rightarrow \cdots$
where each transition map $\Z^i \rightarrow \Z^{i+1}$ is given by
an $(i+1)\times i$ matrix of the form
 $$\begin{pmatrix} n &0 &0 &\cdots &0 &0\\ 1 &0 &0 &\cdots &0 &0\\
 0 &1 &0 &\cdots &0 &0\\ &&&\vdots \\ 0 &0 &0 &\cdots &1 &0 \\ 0 &0 &0
 &\cdots &0 &1 \end{pmatrix},$$
and the class $[B]\in K_0(B)$ corresponds to the
image of $1\in\Z$.

Set $p= y_1x_1\in B$, a proper idempotent. Then set $t_-= x_1\in
U_{-1}$ and $t_+= y_1\in U_1$, so that $t_+t_-= p$ and $t_-t_+
=1$. Hence, the rule $b\mapsto t_+bt_-$ gives a proper corner
isomorphism $\alpha: B\rightarrow pBp$, and Lemma \ref{recognize}
shows that
 $$U= U_{1,n}(k) = \btpmalpha. \qquad\qquad\square$$  }
\end{example}

\begin{example}\label{Uinfty} {\rm
 Let $k$ be a field, and note that there are natural inclusions
$$U_{1,1}(k) \subset U_{1,2}(k) \subset U_{1,3}(k) \subset
\cdots$$
 among the algebras $U_{1,n}(k)$. Set $U_\infty(k)=
 \bigcup_{n=1}^\infty U_{1,n}(k)$, which is a
 simple algebra (e.g., \cite[Theorem 4.3]{agp}). We may also view
 $U_\infty(k)$ as the $k$-algebra with an infinite sequence of
 generators $x_1,y_1,x_2,y_2,\dots$ and relations $x_iy_j=
 \delta_{i,j}$ for all $i,j$. This algebra is $\Z$-graded as
 before, with the $x_i$ having degree $-1$ and the $y_i$ degree
 $1$. Set $B= U_\infty(k)_0$, which is the $k$-linear span
 of the set
$$\{ y_\mu x_\nu \mid \mu,\nu\in \{1,\dots,n\}^\ell,\ n\in\N,\
\ell\in\Zplus \}.$$
 In the present case, $B$ is an ultramatricial
$k$-algebra isomorphic to a direct limit of the algebras
 $$M_{n^n}(k)\times M_{n^{n-1}}(k)\times \cdots\times M_n(k)\times k.$$
Here $K_0(B)$ is isomorphic to the direct limit of a sequence
$\Z^2 \rightarrow \Z^3 \rightarrow \Z^4 \rightarrow \cdots$ with
transition maps
 $$\begin{pmatrix} n &n &n^2 &n^3 &\cdots &n^{n-2} &n^{n-1} \\
 1 &1 &n &n^2 &\cdots &n^{n-3} &n^{n-2} \\
 0 &1 &1 &n &\cdots &n^{n-4} &n^{n-3} \\
 &&&&\vdots \\
 0 &0 &0 &0 &\cdots &1 &1 \\
 0 &0 &0 &0 &\cdots &0 &1 \end{pmatrix} ,$$
and $[B]$ corresponds to $\left( \begin{smallmatrix} 1\\1
\end{smallmatrix} \right) \in \Z^2$. If we define $p$, $t_\pm$,
 $\alpha$ exactly as in Example \ref{U1n}, we conclude from Lemma
 \ref{recognize} that
$$U_\infty(k) = \btpmalpha. \qquad\qquad\square$$  }
\end{example}

\section{Fractional skew monoid rings versus corners of skew group
rings}\label{skewmonfromskewgrp}

Paschke \cite{Paschke} and R\o rdam \cite[Section
2]{rordamclassif} have shown that a C*-algebra crossed product by
an endomorphism corresponds naturally to a corner in a crossed
product by an automorphism. In other words, the C*-algebra
versions of fractional skew monoid rings $\zpazp$ are isomorphic
to corners $e(B*_{\alpha'} \Z)e$ in certain skew group rings. This
leads us to ask whether, in general, our rings $\sopas$ should
appear as corner rings $e(B*G)e$, where $B*G$ is some skew
group ring over the group $G=S^{-1}S$. This is always the case
when $G$ is abelian, as we prove in Proposition
\ref{genunit}. We prepare the way by studying corner rings of the
form $e(A*G)e$ (for $G=S^{-1}S$ as above), and showing that they
fall into the class of fractional skew monoid rings under
appropriate conditions on the action.

\begin{noname}\label{data3} {\rm
Let $A$ be a unital ring, $G$ a group, and $\alpha
:G\rightarrow \aut(A)$ an action. Assume that $S$ is
a submonoid of $G$ with $G=S^{-1}S$, and let $R=A\ast _{\alpha
}G$. Suppose that there exists a nontrivial idempotent $e\in A$
such that $\alpha _s(e)\leq e$ for all $s\in S$.  }
\end{noname}

\begin{lemma}
\label{tech1} Under the above assumptions, the following hold:

{\rm (a)} The action $\alpha$ restricts to an action $\alpha
':S\rightarrow \ndr(eAe)$ by corner
isomorphisms.

{\rm (b)} There are natural monoid morphisms $\sop \rightarrow
eRe$, given by $s\mapsto es^{-1}$, and
$S\rightarrow eRe$, given by $t\mapsto te$, satisfying the
conditions {\rm (1)--(4)} in Definition {\rm \ref{newring}} with
respect to
$\alpha '$ and the inclusion map $\phi :eAe\rightarrow
eRe$.
\end{lemma}

\begin{proof}
(a) This is clear from the hypothesis on $e$.

(b) Notice that, since $e\le \alpha _s^{-1}(e)$ for all $s\in S$,
we have
$es^{-1}=es^{-1}\alpha _s(e)\in eRe$ and
$(es^{-1})(et^{-1})=e(ts)^{-1}$ for $s,t\in S$. Similarly, $se\in
eRe$ and
$(se)(te)=(st)e$. So, the defined maps are monoid morphisms. It is
straightforward to check conditions (1)--(4) in Definition
\ref{newring}.
\end{proof}

Because of Lemma \ref{tech1}, we have the data to construct a
fractional skew monoid ring of the form $\sopasprime$. Since the
maps $\alpha'_s= \alpha_s|_{eAe}$ are injective for all $s\in S$,
the ring homomorphism $eAe\rightarrow \sopasprime$ going with the
construction of $\sopasprime$ is injective by Corollary
\ref{kernel}. Hence, we identify $eAe$
with its image in
$\sopasprime$, as in (\ref{alphasinject}).

\begin{proposition}\label{isomap}
Under the assumptions of {\rm (\ref{data3})}, the rings
$\sopasprime$ and
$e(A\ast _{\alpha }G)e$ are isomorphic as $G$-graded rings.
\end{proposition}

\begin{proof}
By the universal property of $\sopasprime$, there exists a unique
ring homomorphism $\psi :\sopasprime\rightarrow e(A\ast _{\alpha
}G)e$ such that $\psi (s_{-}at_+)=(es^{-1})a(te)$ for all $s,t\in
S$ and $a\in eAe$. Clearly,
$\psi$ is $G$-graded. To see that $\psi$ is onto, consider
$e(ag)e\in e(A*G)e$ where $a\in A$ and $g\in G$, and write
$g=s^{-1}t$ for some $s,t\in S$. Then we have
$$ e(ag)e=eas^{-1}te=(es^{-1})(\alpha_s(ea)\alpha _t(e))(te)\in
\psi (\sopasprime),$$
which proves that $\psi$ is onto. It only remains to check that
$\psi$ is one-to one.

Since $\psi$ is $G$-graded, we only have to check that $\psi
(s_-at_+)=0$ implies $a=0$, when $s,t\in S$ and $a\in
p_s(eAe)p_t$. Note that $p_s= \alpha'_s(1_{eAe})= \alpha_s(e)$,
and likewise $p_t= \alpha_t(e)$, so that $a= \alpha
_s(e)a\alpha _t(e) $.  Now
$$0=(es^{-1})a(te)=e\alpha _s^{-1}(a\alpha _t(e))(s^{-1}t)=\alpha
_s^{-1}(\alpha _s(e)a\alpha _t(e))(s^{-1}t)=\alpha
_s^{-1}(a)(s^{-1}t),$$
whence $\alpha _s^{-1}(a)=0$ and $a=0$, as
desired.
\end{proof}

The following procedure gives a generic way to obtain a situation
as in (\ref{data3}).

\begin{example}
\label{commex} Let $\alpha \colon G\rightarrow \aut (A)$ be an
action of an abelian group $G$ on a unital ring $A$, and let $e$ be
an idempotent in $A$. Set $S :=\{ s\in G\mid \alpha _s(e)\le e\}$.
Then $S$ is a submonoid of $G$ and $G' := S^{-1}S$ is a subgroup of
$G$ acting on $A$ via $\alpha$. Moreover,
$e(A*_{\alpha}G')e\cong \sopasprime$, where $\alpha '\colon
S\rightarrow \ndr(eAe)$ is the induced action of $S$ on
$eAe$ by corner isomorphisms.
\end{example}

\begin{proof}
It is clear that $S$ is a submonoid of $G$, and we can apply
Proposition \ref{isomap} to get the result.
\end{proof}

Now we go in the reverse direction, looking for a representation
of a fractional skew monoid ring $\sopas$ as a corner ring of a
skew group ring. For this to hold, we shall need $S$ to be
abelian. We follow ideas of R\o rdam \cite{rordamclassif}.

\begin{noname}\label{data4} {\rm
Let $A$ be a unital ring, $G$ a group and $S$ a
submonoid of $G$ such that $G=S^{-1}S$. Let $\alpha :S\rightarrow
\ndr(A)$ be an action of $S$ on $A$ by corner
isomorphisms, and for $s\in S$ let $p_s$ denote the idempotent
$\alpha _s(1)$.

We consider the pre-order on $S$ given by $s\le t$ if and only if
there exists $s'\in S$ such that $t=s's$, that is, if and only if
$ts^{-1}\in S$. With this pre-order,
$S$ can be considered as a small category, with at most one arrow
between two given objects. Note that since $G=S^{-1}S$, we have
the left Ore condition, and so $S$ is upward directed.

We define a functor $F$ from the category $S$ to the category of
non-unital rings as follows. For each $s\in S$, set $F(s)=A$. If
$s\le t$ in $S$, then $t=s's$ for some (unique) $s'\in S$; we set
$f_{t,s}= \alpha_{s'}$, viewed as a morphism $A= F(s)\rightarrow
A= F(t)$, and define $F(s{\rightarrow}t)= f_{t,s}$. Let $\A$ be
the colimit of
$F$. Since $S$ is directed, this
colimit is in fact a direct limit, and since all the maps
$f_{t,s}$ are injective, the colimit maps $\rho _s \colon
F(s)\rightarrow
\A$ are injective for all $s\in S$. On the other hand, the maps
$f_{t,s}$ are typically non-unital, and so $\A$ may well be a
non-unital ring.  }
\end{noname}

\begin{lemma}\label{aprounit}
Under the above assumptions, set $q_s=\rho _s(1)$ for each $s\in
S$. Then the following properties hold:

{\rm (a)} If $s,t\in S$ and $s\le t$ then $q_s\le q_t$.

{\rm (b)} $\A=\bigcup _{s\in S} q_s\A q_s$.

{\rm (c)} For each $s\in S$, the map $\rho_s$ gives an isomorphism
$A \overset{\cong}\longrightarrow q_s\A q_s$.
\end{lemma}

\begin{proof}
(a) We have $t=s's$ for some $s'\in S$. Then $q_{s}=\rho
_s(1)=(\rho _t\circ f_{t,s})(1)=\rho _t(\alpha _{s'}(1))\le \rho
_t(1)=q_t$.

(b) Since $\A$ is the direct limit of the family
$\{F(s)\}_{s\in S}$, we have $\A=\bigcup _{s\in S}\rho
_s(F(s))$. Each $\rho_s(F(s))= \rho_s(1A1) \subseteq q_s\A q_s$,
and (b) follows.

(c) We have already observed that $\rho_s(A) \subseteq q_s\A q_s$.
To establish the reverse inequality, it suffices to show that
$q_s\rho_t(A)q_s\subseteq \rho_s(A)$ for each $t\in S$. Given $t$,
we may choose $u\in S$ with $s,t\le u$, since $S$ is upward
directed. Then $\rho_t(A)= \rho_uf_{u,t}(A)\subseteq \rho_u(A)$.
Moreover, $u=s's$ for some $s'\in S$, and $\alpha_{s'}(A)=
\alpha_{s'}(1)A\alpha_{s'}(1)$ because $\alpha_{s'}$ is a corner
isomorphism. Thus,
\begin{align}
\nonumber \rho_s(A) &= \rho_uf_{u,s}(A)= \rho_u\alpha_{s'}(A)=
\rho_u \bigl( \alpha_{s'}(1)A\alpha_{s'}(1) \bigr) \\
\nonumber &= \rho_uf_{u,s}(1)\rho_u(A)\rho_uf_{u,s}(1)=
q_s\rho_u(A)q_s \supseteq q_s\rho_t(A)q_s.
\end{align}
This shows that $\rho_s(A)= q_s\A q_s$. Since $\rho_s$ is
injective, (c) is proved.
\end{proof}

\begin{lemma}\label{grpaction}
Under the assumptions {\rm (\ref{data4})}, if $G$ is an abelian
group, then the action
$\alpha$ extends to an action ${\widehat {\alpha }}:G\rightarrow
\aut(\A)$.
\end{lemma}

\begin{proof}
The group $G$ is the universal enveloping group of $S$, so it
suffices to extend $\alpha$ to a semigroup action $\alsom
\colon S\rightarrow \aut (\A)$.

First, we fix an element $u\in S$, and define $\alsom _u$. We have
ring homomorphisms $\rho_s\circ\alpha_u: F(s)=A \rightarrow \A$
for all $s\in S$. If $s,t\in S$ and $s\le t$, then $t=s's$ for
some $s'\in S$, and
$$(\rho_t\circ\alpha_u)\circ f_{t,s}= \rho_t\circ \alpha_u\circ
\alpha_{s'}= \rho_t\circ
\alpha_{s'}\circ \alpha_u= \rho_t\circ f_{t,s}\circ \alpha_u=
\rho_s\circ \alpha_u.$$
 (Here we have used
the commutativity of $S$.) By the universal property of the
colimit, there exists a unique ring endomorphism $\alsom _u\colon
\A\rightarrow \A$ such that $\alsom _u\circ\rho_s=
\rho_s\circ\alpha_u$ for all $s\in S$. Since each
$\rho_s\circ\alpha_u$ is injective, so is $\alsom _u$. Moreover,
for $s\in S$ we have
$$\alsom_u\circ\rho_{us}= \rho_{us}\circ\alpha_u= \rho_{us}\circ
f_{us,s}= \rho_s$$
and so $\rho_s(A) \subseteq \alsom_u(\A)$. This shows that $\alsom
_u$ is onto. Now, it is straightforward to show that $\alsom
_{uv}=\alsom _u\circ \alsom _v$ for $u,v\in S$, which completes
the proof.
\end{proof}

\begin{proposition}\label{genunit}
Let $G$ be an abelian group and $S$ a submonoid of $G$ such
that $G=S^{-1}S$. Let $\alpha :S\rightarrow \ndr(A)$ be
an action of $S$ on $A$ by corner isomorphisms. Then there exist a
unital ring $B$, an action $\alsom: G\rightarrow \aut(B)$, and an
idempotent $e$ in $B$ such that $\alsom _s(e)\le e$ for all $s\in
S $ and $\sopas \cong e(B*_{\hat\alpha}G)e$ (as $G$-graded rings).
\end{proposition}

\begin{proof}
Let $\A$ and $\alsom $ be as above. The ring $\A$
will not be unital in general, so let $B$ be the unitization of
$\A$. Extend $\alsom $ to
an action on $B$ in the obvious way, and view $\A $ as a
two-sided ideal of $B$. Let $e=q_1= \rho _1(1)\in B$. It is
obvious then that $\alsom (e)\le e$ for every $s\in S$. Thus,
$\alsom$ restricts to an action $\beta: S\rightarrow \ndr(eBe)$
by corner isomorphisms (Lemma \ref{tech1}), and
$\sop*_{\beta} (eBe)*_{\beta} S \cong
e(B*_{\hat\alpha}G)e$ as $G$-graded rings (Proposition
\ref{isomap}). By Lemma
\ref{aprounit}(c), $\rho_1$ gives a ring isomorphism $\rho:A
\rightarrow q_1\A q_1= eBe$, and since $\alsom _s\circ\rho_1=
\rho_1\circ\alpha_s$ for all $s\in S$, we see that $\rho$
transports the action $\alpha$ to the action $\beta$. Therefore
$\sopas \cong \sop*_{\beta} (eBe)*_{\beta} S$ as $G$-graded rings.
\end{proof}

\section{Simplicity}

We continue the general assumptions of (\ref{data1}) and
(\ref{data2}), and seek conditions on $A$, $S$, and $\alpha$ under
which $R=\sopas$ is a simple ring. In the case of a group action
(i.e., $S=G$ and $\alpha:G\rightarrow \aut(A)$), sufficient
conditions for simplicity are well known \cite[Theorem 2.3]{Mon}:
If $A$ is simple and the action $\alpha$ is outer, then the skew
group ring $A*_\alpha G$ is simple. It turns out that a suitable
modification of the notion of an outer action also leads to
simplicity in our more general situation.

We shall say that a pair $(\alpha_s,\alpha_t)$, where $s,t\in S$,
is \emph{inner} provided there exist elements $u\in p_sAp_t$ and
$v\in p_tAp_s$ such that $uv=p_s$, $vu=p_t$ and $\alpha
_s(x)=u\alpha _t(x)v$ for all $x\in A$. Note that then $\alpha
_s\alpha _t^{-1}(x)=uxv$ for every $x\in p_tAp_t$, and
$\alpha_t\alpha_s^{-1}(x)= vxu$ for all $x\in p_sAp_s$. Let us say
that $\alpha$ is \emph{outer} in case $(\alpha_s,\alpha_t)$ is not
inner for any distinct $s,t\in S$.

We will use the following standard terminology. The \emph{support}
of an element  $r=\sum _x r_x$ in $R=\bigoplus _{x\in G}R_x$ is
the set $\mbox{Supp}(r)=\{x\in G\mid r_x\neq 0$\}. The
\emph{length} of $r$ is the number of elements in the support of
$r$, and is denoted $\length(r)$.

\begin{theorem}\label{simplicity}
If $A$ is simple and $\alpha$ is outer, then $R=\sopas$ is simple.
\end{theorem}

\begin{proof}
Suppose that $R$ is not simple. Let $I$ be a proper nonzero ideal
of $R$, and let $\rho\in I$ be a nonzero element with minimal
length, say length $n$. Write $\rho= \sum_{i=1}^n
(s_i)_-a_i(t_i)_+$ where the $s_i^{-1}t_i$ are distinct elements
of $\sinvs$ and each $a_i$ is a nonzero element of
$p_{s_i}Ap_{t_i}$. Observe that $(s_1)_+\rho(t_1)_-= a_1+
\sum_{i=2}^n \rho_i$ where each $\rho_i$ lies in the
$s_1s_i^{-1}t_it_1^{-1}$-component of $R$. Hence,
$(s_1)_+\rho(t_1)_-= a_1+ \sum_{i=2}^n (u_i)_-b_i(v_i)_+$ where
the $u_i^{-1}v_i$ are distinct elements of $\sinvs$, different
from $1$, and each $b_i\in p_{u_i}Ap_{v_i}$. Moreover, $a_1\ne 0$
implies $(s_1)_+\rho(t_1)_- \ne 0$, and so $(s_1)_+\rho(t_1)_-$
has length $n$ by minimality. Thus, after replacing $\rho$ by
$(s_1)_+\rho(t_1)_-$, we may assume that $s_1=t_1=1$.

Since $A$ is simple, $\sum_{j=1}^m c_ja_1d_j=1$ for some
$c_j,d_j\in A$. Then we can replace $\rho$ by
$$\sum_{j=1}^m c_j\rho d_j= 1+ \sum_{i=2}^n (s_i)_- \bigl(
\sum_{j=1}^m \alpha_{s_i}(c_j)a_i\alpha_{t_i}(d_j) \bigr)
(t_i)_+,$$
 and so we may now assume that $a_1=1$. Of course
$\rho\ne 1$ because $I\ne R$, whence $n\ge 2$. Set $s=s_2$,
$t=t_2$, and $a=a_2\in p_sAp_t$, so that
$$\rho= 1+ s_-at_++ \sum_{i=3}^n (s_i)_-a_i(t_i)_+.$$

For any $x\in A$, we have $x\rho-\rho x\in I$ and
$$x\rho-\rho x= s_-\bigl( \alpha_s(x)a- a\alpha_t(x) \bigr)t_+
+\sum_{i=3}^n (s_i)_-\lozenge_i(t_i)_+$$
 for some elements $\lozenge_i \in p_{s_i}Ap_{t_i}$ that we need
 not specify. Thus $x\rho-\rho x$ has length less than $n$, and so
 $x\rho-\rho x=0$ by the minimality of $n$. Therefore
 $$\alpha_s(x)a= a\alpha_t(x)$$
for all $x\in A$. In particular, $p_sAa= p_sAp_sa= \alpha_s(A)a=
a\alpha_t(A)= aAp_t$.

Since $A$ is simple, $Ap_tA= AaA= A$, and so
$$aAp_s= aAp_tAp_s= p_sAaAp_s= p_sAp_s,$$
 whence there is some $b\in p_tAp_s$ such that $ab=p_s$.
 Similarly, there is some $c\in p_tAp_s$ such that $ca= p_t$. But
 $c=cp_s= cab= p_tb=b$, so that $ba=p_t$. Now
$$a\alpha_t(x)b= \alpha_s(x)ab= \alpha_s(x)p_s= \alpha_s(x)$$
for all $x\in A$, and so we conclude that the pair
$(\alpha_s,\alpha_t)$ is inner. Since $\alpha$ is assumed to be
outer, we must have $s=t$. But then $s_2^{-1}t_2= s^{-1}t= 1=
 s_1^{-1}t_1$, contradicting the distinctness of the
 $s_i^{-1}t_i$. Therefore $R$ is simple.
\end{proof}

\begin{corollary}\label{notequiv}
If $A$ is simple and $p_s \not\sim p_t$ for all distinct $s,t\in
S$, then $R$ is simple. \qquad\qquad$\square$
\end{corollary}

\begin{corollary}\label{df}
If $A$ is a directly finite simple ring, $p\in A$ is a proper
idempotent (i.e., $p\ne 1$), and $\alpha :A\rightarrow pAp$ is a
corner isomorphism, then $\zpazp$ is simple.
\end{corollary}

\begin{proof} The idempotents corresponding to the monoid
homomorphism $\Zplus \rightarrow \ndr(A)$ in this case are the
$\alpha^i(1)$ for $i\in\Zplus$. Since $\alpha(1)=p \ne 1$, we have
$1> \alpha(1)> \alpha^2(1)> \cdots$, and it follows from the
direct finiteness of $A$ that $\alpha^i(1) \not\sim \alpha^j(1)$
for all distinct $i,j\in \Zplus$.
\end{proof}

\section{Purely infinite simplicity}

We recall from \cite{agp} that a simple ring $T$ is said to be
\emph{purely infinite} if every nonzero right ideal of $T$
contains an infinite idempotent. This concept is left-right
symmetric, as the following characterization shows: $T$ is purely
infinite if and only if $(1)$ $T$ is not a division ring; $(2)$
for every nonzero element $a\in T$, there exist elements $x,y\in
T$ such that $xay=1$ \cite[Theorem 1.6]{agp}. For instance, the
Leavitt algebras $V_{1,n}(k)$ and $U_\infty(k)$ are purely
infinite simple rings \cite[Theorems 4.2, 4.3]{agp}. As we have
seen above (Examples \ref{algCK} and \ref{Uinfty}), the
$V_{1,n}(k)$ and $U_\infty(k)$ can be presented in the form
$\zpbzp$. This suggests that fractional skew monoid rings might be
purely infinite simple in some generality. Our goal in this
section is to establish sufficient conditions for a fractional
skew monoid ring $R= \sop*_\alpha A*_\alpha S$ to be a purely
infinite simple ring, under the general assumptions of
(\ref{data1}) and (\ref{data2}).

The following concept will be needed. A ring $T$ is said to be
\emph{strictly unperforated} provided the finitely generated
projective right (or left) $T$-modules enjoy the following
property: If $mA \prec mB$ for some $m\in\N$, then $A\prec B$.
(Here $mA$ denotes the direct sum of $m$ copies of $A$, and the
notation $X\prec Y$ means that $X$ is isomorphic to a
\underbar{proper} direct summand of $Y$.) Stated in terms of
idempotents in matrix rings over $T$, strict unperforation is the
condition $(m{{\cdot}} p\prec m{\cdot} q \implies p\prec q)$, where $m{\cdot} p$
denotes the orthogonal sum of $m$ copies of an idempotent $p$. For
instance, ultramatricial algebras are strictly unperforated
\cite[Theorem 15.24(a)]{vnrr}. Also, any purely infinite simple
ring $T$ is strictly unperforated, because $A\prec B$ for all
nonzero finitely generated projective $T$-modules $A$ and $B$
\cite[Proposition 1.5]{agp}.

\begin{lemma}\label{Paso1}
Assume that $A$ is simple and strictly unperforated, and that
there exists $u\in S$ such that $p_u \ne 1$. For any nonzero
idempotent $e\in A$, there exists $v=u^j\in S$ for some $j\in \N$
such that $p_v\lesssim e$.
\end{lemma}

\begin{proof}
Set $p_i= p_{u^i}= \alpha_u^i(1)$ for $i\ge 0$. Since $A$ is
simple, there exists $m\in\N$ such that $1\prec m{\cdot} e$ and
$1\lesssim m{\cdot} (1-p_1)$. Note that
$$(m+1){\cdot} p_1\lesssim m{\cdot} p_1\oplus 1\lesssim m{\cdot} p_1\oplus
m{\cdot} (1-p_1)\sim m{\cdot} 1.$$
 Applying the isomorphisms $\alpha_u^i: A\rightarrow p_iAp_i$,
 we obtain that $(m+1){\cdot} p_{i+1}\lesssim m{\cdot} p_i$ for all $i$. It
 follows by induction that $(m+1)^i{\cdot} p_i\lesssim m^i{\cdot} 1$ for all
 $i$.

Now choose $j\in\N$ such that $m^{j+1}< (m+1)^j$, and observe that
$$m^{j+1}{\cdot} p_j\prec (m+1)^j{\cdot} p_j\lesssim m^j{\cdot} 1\prec
m^{j+1}{\cdot} e,$$
 whence $m^{j+1}{\cdot} p_j\prec m^{j+1}{\cdot} e$. Therefore $p_j\prec e$,
  because $A$ is strictly
 unperforated.
\end{proof}

The following lemma is a variation on results such as
\cite[Proposition 3.3]{vnrr}.

\begin{lemma}\label{idemgen}
If $T$ is a simple ring containing an idempotent $p\ne 0,1$, then
$T$ is generated (as a ring) by its idempotents.
\end{lemma}

\begin{proof} Let $T'$ be the subring of $T$ generated by the
idempotents. Since $p+pt(1-p)$ is idempotent for any $t\in T$, we
see that $pT(1-p) \subseteq T'$, and likewise $(1-p)Tp \subseteq
T'$. The simplicity of $T$ implies that $T(1-p)T =T$, whence $pTp=
[pT(1-p)][(1-p)Tp] \subseteq T'$, and similarly $(1-p)T(1-p)
\subseteq T'$. Therefore $T'=T$.
\end{proof}

\begin{theorem}\label{spi} Assume that $A$ is a simple, strictly
unperforated ring, in which every nonzero right (left) ideal
contains a nonzero idempotent. Assume also that $\alpha$ is outer,
and that there exists $u\in S$ with $p_u\ne 1$. Then $R=\sopas$ is
a purely infinite simple ring.
\end{theorem}

\begin{proof} The hypothesis that
$p_u\ne 1$ will allow us later to apply Lemma \ref{Paso1}.
Moreover, it implies that $R$ is not a division ring.

 Let $\rho$ be an arbitrary nonzero element of $R$.
Choose $\rho',\rho''\in R$ such that $\rho'\rho\rho''$ is nonzero
and has minimal length for such nonzero products, say length $n$.
Since it suffices to find $x,y\in R$ such that
$x\rho'\rho\rho''y=1$, we may replace $\rho$ by $\rho'\rho\rho''$.
Thus, without loss of generality, all nonzero products
$\sigma\rho\sigma'$ in $R$ have length at least $n$. Now write
$\rho= \sum_{i=1}^n (s_i)_-a_i(t_i)_+$ where the $s_i^{-1}t_i$ are
distinct elements of $S^{-1}S$ and each $a_i$ is a nonzero element
of $p_{s_i}Ap_{t_i}$. As in the proof of Theorem \ref{simplicity},
after replacing $\rho$ by $(s_1)_+\rho(t_1)_-$ we may assume that
$s_1=t_1=1$, so that $\rho=a_1+ \sum_{i=2}^n (s_i)_-a_i(t_i)_+$.

By our hypothesis on idempotents, there exists $a'_1\in A$ such
that $a_1a'_1$ is a nonzero idempotent. By Lemma \ref{Paso1},
there exist $x,y\in A$ such that $xa_1a'_1y=p_v$ for some $v\in
S$. Note that $v_-xa_1a'_1yv_+=1$. Hence, after replacing $\rho$
by $v_-x\rho a'_1yv_+$, we may assume that $a_1=1$. We are thus
done in case $n=1$.

Suppose that $n\ge2$, and set $s=s_2$, $t=t_2$, and $a=a_2\in
p_sAp_t$. Thus,
$$\rho= 1+ s_-at_++ \sum_{i=3}^n (s_i)_-a_i(t_i)_+$$
 at this point. For any idempotent $e\in A$, we have
$$e\rho(1-e)= s_-\alpha_s(e)a \bigl( p_t-\alpha_t(e) \bigr) t_++
 \sum_{i=3}^n (s_i)_-\lozenge_i(t_i)_+.$$
 Since $e\rho(1-e)$ has length less than $n$, it must
 be zero, whence $\alpha_s(e)a(p_t-\alpha_t(e))=0$. Thus,
 $\alpha_s(e)a= \alpha_s(e)a\alpha_t(e)$. A symmetric argument
 involving $(1-e)\rho e$ shows that $a\alpha_t(e)=
 \alpha_s(e)a\alpha_t(e)$, and so $\alpha_s(e)a= a\alpha_t(e)$.

By Lemma \ref{idemgen}, $A$ is generated by its idempotents.
Hence, it follows from the equations $\alpha_s(e)a= a\alpha_t(e)$
that $\alpha_s(x)a= a\alpha_t(x)$ for all $x\in A$. As in the
proof of Theorem \ref{simplicity}, this implies that the pair
$(\alpha_s,\alpha_t)$ is inner, yielding $s=t$ and $s_2^{-1}t_2=
s_1^{-1}t_1$, which contradicts our assumptions. Therefore $n=1$,
and the proof is complete.
\end{proof}

It is perhaps not so surprising that the purely infinite simple
property carries over from $A$ to $R$ under suitable conditions.
More interesting is that $R$ can be purely infinite simple even
when $A$ is directly finite. We single out an important case of
this phenomenon in the following corollary.

\begin{corollary}\label{spicor} Suppose that $A$ is either a
purely infinite simple ring or a simple ultramatricial algebra
over some field. Assume also that $\alpha$ is outer, and that
there exists $u\in S$ with $p_u\ne 1$. Then $R=\sopas$ is a purely
infinite simple ring. \qquad\qquad$\square$
\end{corollary}

\section*{Acknowledgments}

Parts of this work were done during visits of the first author to
the Department of Mathematics of the University of California at
Santa Barbara and to the Departamento de Matem\'at\-i\-cas de la
Universidad de C\'adiz, and during visits of the last three
authors to the Departament de Matem\`atiques de l'Universitat
Aut\`onoma de Barcelona and to the Centre de Recerca Matem\`atica
(UAB). The authors would like to thank these host centers for
their warm hospitality.

\end{document}